\documentclass[12pt]{article}
\usepackage{amsmath}
\usepackage{amssymb}
\usepackage{verbatim}
\newtheorem{thm}{Theorem}

\newcommand{\R}{\mathbb{R}}
\newcommand{\RN}{\mathbb{R}^N}
\newcommand{\dif}[0]{\ensuremath{\,\mathrm{d}}}

\DeclareMathOperator*{\dive}{div}

\begin{document}
\title{ON THE MEAN VALUE PROPERTY FOR THE $p$-LAPLACE EQUATION IN THE PLANE}
\author{Peter Lindqvist \qquad Juan Manfredi}
\date{\footnotesize{Norwegian University of Science and Technology.\; University of Pittsburgh}}
\maketitle

{\small \textsc{Abstract:}\footnote{AMS classification 35J92,
35J62.} \textsf{We study the $p$-Laplace equation in the plane and prove that the mean value property holds directly for the solutions themselves. This removes the need to interpret the formula in the viscosity sense via test functions. The method is based on the hodograph representation.}

\section{Introduction}

Harmonic functions can be characterized by the \lq\lq asymptotic\rq\rq\  mean value property
$$ u(x)\; =\; \frac{1}{|B(x,\varepsilon)|}\int_{B(x,\varepsilon)}\!u(y)\,\dif y + o(\varepsilon^2)\qquad \text{as}\qquad \varepsilon \to 0$$
valid at each point $x$ in the domain of definition.
The expansion is an equality  for harmonic functions with $o(\varepsilon^2) = 0$. It is an  exercise to prove that only the asymptotic expansion is needed to conclude that a $C^{2}$-function is harmonic, and it is an very interesting exercise to prove it for continuous functions.\par
 In [MRP] it was proved that a similar property holds for the solutions of the $p$-Laplace equation 
$$\dive(|\nabla u|^{p-2}\nabla u)\, = \, 0,$$
which is the Euler-Lagrange equation for the variational integral
$$ I(u) = \int_{\Omega} | \nabla u|^{p}\,dx,$$
where $1<p<\infty.$
Indeed, a function $u \in C(\Omega)\cap W^{1,p}_{loc}(\Omega)$ is a solution in the domain $\Omega$ in $\RN$ of the $p$-Laplace equation if and only if the asymptotic expansion

\begin{equation}
\label{meann}
u(x) = \tfrac{p-2}{p+n}\; \left(\frac{\underset{\overline{B}(x,\varepsilon)}{\max}\,u + \underset{\overline{B}(x,\varepsilon)}{\min}\,u}{2}\right)\\ \!\!\!\
+\tfrac{2+n}{p+n}\;\left(\frac{1}{|B(x,\varepsilon)|}\int_{B(x,\varepsilon)}u(y)\,\dif y \right)\! + o(\varepsilon^2)
\end{equation}
holds at each point $x \in \Omega$, as $\varepsilon \to 0.$ However, the expansion was proved to be valid in the \emph{viscosity sense}, which means that, strictly speaking, $u$ has to be replaced by test  functions on which the pointwise calculations are performed. In the extreme case $p = \infty$, when the equation reads
$$\sum_{i,j=1}^{n}\frac{\partial u}{\partial x_{i}} \frac{\partial u}{\partial x_{j}}
\frac{\partial^{2} u}{\partial x_{i}\partial x_{j}}\;=\; 0,$$
a counter example shows that the expected  formula
$$ u(x) =  \frac{\underset{\overline{B}(x,\varepsilon)}{\max}\,u + \underset{\overline{B}(x,\varepsilon)}{\min}\,u}{2} + o(\varepsilon^2)$$
cannot be interpreted literally as it stays; indeed, test functions are needed. (The Aronsson function $u(x,y) = x^{4/3} - y^{4/3}$ will do as an example.) For finite values of $p$ we do not know about any counter example.

The first derivatives are known to be continuous, but the lack of second derivatives is crucial. It follows essentially from the Taylor expansion that if $u \in C^{2}(\Omega)$ and $\nabla u(x) \not = 0$ then (\ref{meann}) holds directly for $u$ itself, i.e., without any viscosity interpretation. (In fact, it holds also when $\nabla u(x) = 0$ if, in addition, $D^2u(x) = 0.$) Since solutions of the $p-$Laplace equation are even real analytic at points where $\nabla u \not= 0$ according to [L, p. 208], the whole problem is at the critical points. Unfortunately, in space nothing seems to be known about the critical points.

In the plane, much more structure is available due to the fact that the (complex) gradient $u_{x} - i u_{y}$ is a quasiregular mapping. This fundamental result of Bojarski-Iwaniec implies that \emph{the critical points are isolated}, unless $u$ is constant, cf. [BI]. We shall use the hodograph representation given in [IM] to prove that in the plane the mean value formula holds pointwise at least in the range $$1<p < p_0 = 9.52520797...$$

\begin{thm}
Suppose that $\Omega$ is a domain in the plane $\R^{2}$ and let $1< p < p_{0} = 9.52...$
A function $u\in C(\Omega)\cap W_{loc}^{1,p}(\Omega)$ is a solution to the $p$-Laplace equation if and only if the expansion
\begin{equation*}
u(x) = \frac{p-2}{p+2}\; \left(\frac{\underset{\overline{B}(x,\varepsilon)}{\max}\,u + \underset{\overline{B}(x,\varepsilon)}{\min}\,u}{2}\right)\\
+\frac{4}{p+2}\;\left(\frac{1}{\pi \varepsilon^{2}}\int_{B(x,\varepsilon)}u(y)\,dy \right) + o(\varepsilon^2)
\end{equation*}
holds at each point $x \in \Omega$ as $\varepsilon \to 0.$
\end{thm}

The number $p_{0}$ restricting the range of validity seems to be an artifact; it is the root of an auxiliary equation. To get beyond $p_0$ would require far more complicated calculations. We  recall that $u \in C(\Omega)\cap W_{loc}^{1,p}(\Omega)$ is a solution of the $p$-Laplace equation if and only if

$$\int_{\Omega}\langle |\nabla u|^{p-2}\nabla u,\nabla \varphi \rangle \dif x  = 0$$
for all  $\varphi \in C_{0}^{\infty}(\Omega).$ In fact, $u \in C_{loc}^{1,\alpha}(\Omega)$ and, in the plane, a complete regularity characterization is given in [IM].

Our method of proof relies on a separation of the first (and worst) term in the hodograph representation. This term can, with some care, be handled without destruction of the remainder.

\paragraph{Acknowledgement}{\footnotesize This research was done when one author visited the University of Pittsburgh. He  thanks the Department of Mathematics.}

\section{Preliminaries}

We sketch the hodograph method, for which we refer the reader to consult [IM]. It is based on the fact that if $u = u(x,y)$ is a solution to the $p$-Laplace equation in the plane domain $\Omega$, then the complex gradient
$$f(z) = \frac{\partial u}{\partial z} = \frac{1}{2}\left(\frac{\partial u}{\partial x} - i  \frac{\partial u}{\partial y}\right)$$
is a \emph{quasiregular mapping} according to [BI]. (By abuse of notation, $f(z)$ is used for $f(z,\overline{z}).$) It follows that \emph{the zeros of $f$ are isolated}, unless $f$ is a constant. Assume that $z_0 = x_0 + i y_0$ is a critical point, i.e., $f(z_0) = 0.$ By adding a constant, we may assume that $u(x_0,y_0) = 0.$ For some integer $n \geq 1$ we have the Stoilow representation
$$f(z) = \chi(z)^n$$
in a small neighbourhood $|z-z_0| < \varrho,$ where $\chi$ is a quasiconformal mapping, i.e., it is also injective.

We shall reproduce a formula from [IM]. Let us for simplicity take $z = 0$ so that now $f(0) = 0,\, \chi(0) = 0.$ The inverse mapping defined by
$$ \zeta = \chi(z) \quad  \Longleftrightarrow \quad z = H(\zeta)$$
in a suitable neighbourhood has the representation formula
\begin{align*}
H(\zeta)\;&= \sum_{k=n+1}^{\infty}\left(A_k\zeta^k+\varepsilon_k\overline{A}_k\overline{ \zeta}^k\right)\zeta^{-n}|\zeta|^{\lambda_{k}+n-k}\\
&=\overbrace{\left\{A_{n+1}\frac{\zeta}{|\zeta|}+\varepsilon_{n+1}\overline{A}_{n+1}\left(\frac{\overline{\zeta}}{\zeta}\right)^{n}\right\}|\zeta|^{\lambda_{n+1}}}^{A(\zeta)}\\
&+\left\{A_{n+2}\left(\frac{\zeta}{|\zeta|}\right)^{n+2}+\varepsilon_{n+2}\overline{A}_{n+2}\left(\frac{\overline{\zeta}}{\zeta}\right)^{n+2}\left(\frac{|\zeta|}{\zeta}\right)^{n}\right\}|\zeta|^{\lambda_{n+2}}\\
&+\cdots = A(\zeta) + R(\zeta),
\end{align*}
where the parameters are given by
\begin{align*}
\lambda_k = \lambda_k^{(n)} &= \frac{-np+ \sqrt{4k^2(p-1)+(p-2)^2}}{2}\\
\varepsilon_k = \varepsilon_k^{(n)} &= \frac{\lambda_k + n - k}{\lambda_k + n + k}
\end{align*}
for $k = n+1,n+2,\cdots.$ It is important that $|\varepsilon_k| < 1.$ Thus, if $A_k \not= 0,$ the $k^{th}$ term can have no other zeros than $\zeta =0.$ We record that $0 < \varepsilon_3 <1,$ when $p >2.$ If
\begin{equation}
\label{conv}
\sum_{k=n+1}^{\infty} k|A_{k}|^2 \: < \; \infty
\end{equation}
this formula produces, according to Theorem 2 in [IM], all solutions of the $p$-Laplace equation in a sufficiently small neighbourhood of the critical point $0$.

\subsection{About the cases $n = 1,2,3,\cdots$ for $f(z) = \chi(z)^n$}

A consequence for the second derivatives is that
$$\sum |D^2u(z)|\; \leq \; B_1|z|^{\frac{n}{\gamma_{n}} - 1}$$
in a small neighbourhood of $0$, where
$$\frac{\gamma_n}{n} = \frac{1}{2}\biggl(\sqrt{4(1+\frac{1}{n})^2(p-1)+(p-2)^2}-p\biggr).$$
It follows that $u$ has H\"{o}lder continuous second derivatives, if
$$\frac{n}{\gamma_{n}} > 1.$$
A calculation shows that this holds when
\begin{align*}
n&=1  &  1< p&<2\\
n&=2 &  1< p&<9\\
n&\geq 3  &  1 < p&< \infty.
\end{align*}
However, for our purpose it is impossible to know $n$ in advance. Therefore the constant $p_0$ in our theorem is determined from the most difficult case $n=1.$ 

\paragraph{The case $n=1,\quad f(z) = \chi(z)$}

Let us keep $p   > 2.$ Now we have $A_2 \not = 0,\, \zeta = f(z), \, z = H(\zeta)$ and
\begin{gather}
z = H(\zeta) \quad = \quad \overbrace{\left\{A_{2}\frac{\zeta}{|\zeta|}+\varepsilon_{2}\frac{\overline{\zeta}^2}{|\zeta|\zeta}\right\}|\zeta|^{\lambda_{2}}}^{A(\zeta)}\\
+\left\{\left(A_{3}\left(\frac{\zeta}{|\zeta|}\right)^{3}+\varepsilon_{3}\overline{A}_{3}\left(\frac{\overline{\zeta}}{\zeta}\right)^{3}\right)\frac{|\zeta|}{\zeta}\right\}|\zeta|^{\lambda_{3}} + \cdots
= A(\zeta) + R(\zeta),\nonumber
\end{gather}
where the exponents
$$\lambda_k = \frac{-p + \sqrt{4k^2(p-1)+(p-2)^2}}{2}$$
increase with $k = 2,3,\cdots.$ By assumption $A_2 \not = 0.$ If all the other $A_k$'s are $0$, we are done. The reason is the symmetry $z = A(\zeta) = -A(-\zeta),$
which implies that, upon inversion, $a(z) = -a(-z).$ We have used the notation
$$ z = A(\zeta)\quad \Longleftrightarrow \quad \zeta = a(z)$$
for this special function, which will be used below. Now it is easy to verify that the corresponding $p$-harmonic function $\mathfrak{A}$ \; (i) has mean value 0 and (ii) $\max\mathfrak{A}  = -\min \mathfrak{A}$, taken over a disc $B(0,\varepsilon)$, so that $\max  \mathfrak{A}+ \min  \mathfrak{A}= 0.$ This function appears in [A].
It follows that eqn (\ref{meann}) holds in this case, even with $o(\varepsilon^2) = 0.$ Our method is based on this fact.

If all the remaining $A_k$ are not $0$, there is a smallest $k\geq3$ for which $A_k \not = 0.$ The worst case is $A_3 \not = 0,$ which we now consider. Let us write the remainder as
\begin{align*}
R(\zeta) &= |\zeta|^{\lambda_3}\Bigg\{\left(A_{3}\left(\frac{\zeta}{|\zeta|}\right)^{3}+\varepsilon_{3}\overline{A}_{3}\left(\frac{\overline{\zeta}}{\zeta}\right)^{3}\right)\frac{|\zeta|}{\zeta} \\
&\qquad  \qquad + \sum_{k=4}^{\infty}\left(A_k\left(\frac{\zeta}{|\zeta|}\right)^k+ \varepsilon_k\overline{A}_k\left(\frac{\overline{\zeta}}{\zeta}\right)^k\right)\frac{|\zeta|}{\zeta}\;|\zeta|^{\lambda_{k}-\lambda_{3}} \Bigg\} 
\end{align*}
The first term in braces is dominating for small values of $|\zeta|$, because the series converges and the powers $\lambda_{k}-\lambda_{3}$ are positive. More precisely,
$$(1-\varepsilon_3)|A_3|\quad \leq \quad \Biggl\lvert\left(A_{3}\left(\frac{\zeta}{|\zeta|}\right)^{3}+\varepsilon_{3}\overline{A}_{3}\left(\frac{\overline{\zeta}}{\zeta}\right)^{3}\right)\frac{|\zeta|}{\zeta} \Biggr\rvert .$$
For the series Cauchy's inequality yields the bound
$$\lvert \sum_{k=4}^{\infty}\cdots \rvert ^{2}\; \leq\; 16\sum_{k=4}^{\infty}|A_k|^2\sum_{k=4}^{\infty}|\zeta|^{2(\lambda_{k}-\lambda_{3})},$$
and the exponents in the majorant series are positive and $\asymp 2k\sqrt{p-1}$ as $k \to \infty.$ Recall also (\ref{conv}).  Thus the sum converges and its limit is zero as $\zeta \to 0.$  It follows that  $R(\zeta)$ is real analytic and zerofree in a small punctured disc $0 < |\zeta| < \delta.$

\subsection{The perturbation of the main term, $n = 1$}

So far, we have the setup
\begin{equation*}
\begin{cases}z\, =
A(\zeta) + R(\zeta)\quad \Longleftrightarrow \quad \zeta = f(z)\\
w = A(\xi) \qquad \qquad \, \Longleftrightarrow \quad \xi = a(w).
\end{cases}
\end{equation*}
Let $\zeta = f(z).$  Then
$$ z- R(\zeta) = A(\zeta) \quad \Longleftrightarrow \quad \zeta = a(z-R(\zeta)).$$
It follows that 
$$\zeta = f(z) = a(z) + [a(z-R(\zeta))-a(z)]_{\zeta=f(z)}.$$
We claim that
$$ \zeta = f(z) = a(z) + O(|z|^{\frac{\lambda_3}{\lambda_2\lambda_2}}).$$
Inded, both $a(z)$ and $f(z)$ have the H\"{o}lder exponent $1/ \lambda_2$. Hence
the perturbation term above is
$$|a(z\!-\!R(\zeta))-a(z)| \leq C_2|R(\zeta)|^{\frac{1}{\lambda_2}} \leq C_2C_3|\zeta|^{\frac{\lambda_3}{\lambda_2}} = C_2C_3|f(z)|^{\frac{\lambda_3}{\lambda_2}}
\leq C_2C_3C_2'|z|^{\frac{\lambda_3}{\lambda_2\lambda_2}},$$
and the claim follows.

From this we can further conclude that 
$$\boxed{u(x,y) = \mathfrak{A}(x,y) + O(r^{1+\frac{\lambda_3}{\lambda_2\lambda_2}})}$$
for $r = \sqrt{x^2+y^2}$ small enough. Notice that if 
\begin{equation}
\label{lambda}
\frac{\lambda_3}{\lambda_2\lambda_2}\;> \; 1
\end{equation}
then the error term is of the order  $$O(r^{2+\alpha})$$ for some $\alpha > 0$. It is  inequality (\ref{lambda}) that yields our exponent $p_0 = 9.5...$ and it is under this condition that we can prove the mean value formula. Explicitely, the inequality reads
$$\frac{-p+\sqrt{36(p-1)+(p-2)^2}}{2}\; > \; \frac{\bigl(-p+\sqrt{16(p-1)+(p-2)^2}\bigr)^2}{4}$$
and this holds in the range $1 < p < p_0 = 9.52....$ Upon some manipulations, the number $p_0$ appears as a root of an algebraic equation of the $6^{th}$ degree.

\section{Verification of the Mean Value Formula}

Continuing the case $n=1$, we start from 
$$u(x,y) = \mathfrak{A}(x,y) +\mathfrak{e}(x,y)\quad\text{where}\quad |\mathfrak{e}(x,y)| \leq Cr^{2+\alpha}.$$
For the mean value we have
\begin{gather*}
\frac{1}{\pi \varepsilon^2}\iint_{B(0,\varepsilon)}\!u \, \dif x\dif y\;=\;
\frac{1}{\pi \varepsilon^2}\iint_{B(0,\varepsilon)}\!(u -\mathfrak{A})\, \dif x\dif y\\ \;=\;\frac{1}{\pi \varepsilon^2}\iint_{B(0,\varepsilon)}\!\mathfrak{e}\, \dif x\dif y \;=\; O(\varepsilon^{2+\alpha}),
\end{gather*}
since the mean value of $\mathfrak{A}$ is zero by symmetry.
Using the symmetry again,  we can estimate
\begin{gather*}
\max u + \min u\; = \;\max(\mathfrak{A} + \mathfrak{e}) + \min(\mathfrak{A} + \mathfrak{e}) \\ \leq\; \max\mathfrak{A} +\max \mathfrak{e} + \min\mathfrak{A} + \max\mathfrak{e} \:=\; 
 2\max \mathfrak{e}\; \leq C\varepsilon^{2+\alpha}
\end{gather*}
and from below we obtain in the same way
$$\max u + \min u \;\geq \;- C\varepsilon^{2+\alpha}.$$
In conclusion, the mean value of $u$ and $\tfrac{\max u + \min u}{2}$ are both of order $o(\varepsilon^2)$ and $u(0,0) = 0$. Therefore the formula in Theorem 1 is valid\footnote{We assumed that $A_3 \not = 0$, but since $\lambda_k > \lambda_3$ when $k>3$, inequality (\ref{lambda}) certainly holds.}.

This was the case $n=1.$ In the case $n=2$ we already have the result for $1 < p < 9$. The same method improves the bound $9$ even to a number $> p_0.$ Now the relevant parameters are 
$$\lambda_k^{(2)} = \frac{-2p+ \sqrt{4k^2(p-1)+(p-2)^2}}{2},\qquad k =3,4,\cdots.$$
A pretty similar calculation\footnote{Now $f(z) = \chi(z)^2,$ but the square will cancel in the calculations.} leads to the inequality
$$\frac{\lambda^{(2)}_4}{\lambda^{(2)}_3\lambda^{(2)}_3}\;> \; 1$$
in the place of inequality (\ref{lambda}). We omit the details. Finally, the cases $n \geq 3$ are already clear.---This concludes the verification of the mean value formula.

\end{document}